 \documentstyle[graphicx,amscd,amssymb,verbatim,12pt,righttag,color]{amsart}
 \newtheorem{theorem}{Theorem}[section]
 
 \newtheorem{Prop}[theorem]{Proposition}
 \newtheorem{Lem}[theorem]{Lemma}
 \newtheorem{Cor}[theorem]{Corollary}
 
 \newtheorem{Example}[theorem]{Example}

 \numberwithin{equation}{section}
 
\begin{document}

%
%\date{}
\author{Chun-Kit Lai}
\address{Department of Mathematics, The Chinese University of Hong Kong,  Hong Kong}
\email{cklai@@math.cuhk.edu.hk}
%\thanks {}

\title  {On Fourier frame of absolutely continuous measures}
\date{\today}
\keywords {Absolutely continuity, Bernoulli convolutions, Beurling
densities, Fourier frames, self-similar measures.}
\subjclass{Primary 42C15, 28A25;
 Secondary 28A80.}
\thanks{ The research is supported in part by the HKRGC Grant.}
\maketitle

\begin{abstract}
Let $\mu$ be a compactly supported absolutely continuous probability
measure on ${\Bbb R}^n$, we show that $L^2(K,d\mu)$ admits a Fourier
frame  if and only if its Radon-Nikodym derivative is bounded above
and below almost everywhere on the support $K$. As a consequence, we
prove that if $\mu$ is an equal weight absolutely continuous
self-similar measure on ${\Bbb R}^1$  and $L^2(K,d\mu)$ admits a
Fourier frame, then the density of $\mu$ must be a characteristic
function of self-similar tile. In particular, this shows for almost
everywhere $1/2<\lambda<1$, the $L^2$ space of  the
$\lambda$-Bernoulli convolutions cannot admit a Fourier frame.
   \end{abstract}

\begin{section}{\bf Introduction}

 Let $H$ be a Hilbert space, a sequence of vectors $\{e_i\}_{i\in{\Bbb Z}}$ is called a \textit{frame} if there exists $A,B>0$ such that for any $f\in H$,
\begin{equation} \label{Fourier frame}
 A\|f\|^2\leq \sum_{i\in{\Bbb Z}}|\langle f, e_i\rangle|^2\leq B\|f\|^2.
\end{equation}
The constants $A$ and $B$ are called the the \textit{lower frame bound} and \textit{upper frame bound} respectively. Frame is a natural generalization of orthonormal basis (where $A = B= 1$). It is easily seen from the lower bound that frame is complete in $H$. If  $\{e_i\}_{i\in{\Bbb Z}}$ only satisfies the upper bound in (\ref{Fourier frame}), we call $\{e_i\}_{i\in{\Bbb Z}}$ a \textit{Bessel's sequence}.

\bigskip

The study of frames on Hilbert space was first introduced by Duffin and Schaeffer [DS] in connection with the non-harmonic Fourier series, and has received a lot of attention. Nowadays,  frames are regarded as ``overcomplete bases" since they provide basis-like (though non-unique) expansion of vectors. Because of its redundancy, it provides better stability compared to orthonormal basis. For the Hilbert space of the $L^2$ space of functions,  various kinds of frames such as Fourier frames, Gabor frames, and wavelet frames have been studied. They have close links with time-frequency analysis, sampling theory, and wavelets. One may refer to [Chr] and [G] for some excellent expositions.

\bigskip

 In this paper, we will focus on the Fourier frame. Let $\mu$ be a compactly supported probability measure on ${\Bbb R}^n$. As we will deal with $L^2$ space of different measures on different supports,  we use $L^2(K,d\mu)$ to denote the $L^2$ space of the measure $\mu$ with supp$\mu=K$. In particular, $L^2(K,dx)$ is the $L^2$ space of the Lebesgue measure supported on $K$.  We say that a sequence of  complex exponentials $\{e^{2 \pi i \langle\lambda,\cdot\rangle}\}_{\lambda\in\Lambda}$ is a \textit{Fourier frame} of $L^2(K,d\mu)$ (or just $\mu$) if it is a frame on the Hilbert space $L^2(K,d\mu)$ and  $\Lambda$ is called a \textit{frame spectrum}.

\medskip

Traditionally, the studies of the Fourier frame focus on the case where $\mu$ is the Lebesgue measure supported on $[0,1]$. The work of Landau, Jaffard, and Seip ([Lan], [Ja], [S]) relates the frame spectrum of $L^2([0,1],dx)$ closely with the Beurling densities (see Section 2).  Ortega-Cerd\`{a} and Seip recently completely characterized the frame spectrum $L^2([0,1],dx)$ using de Branges' theory of Hilbert space of entire functions [OS].

\medskip

The more recent study has discovered that some other probability measures can also admit exponential orthonormal bases. One of the surprising results is by Jorgensen and Pedersen [JP], they discovered that the Cantor measures with even contractions ratio admit an exponential orthonormal basis, while those with odd contractions do not. It is still open whether the one-third Cantor measure will admit any Fourier frame ([DHSW], [DHW]).

\bigskip

In the following, we study the existence of Fourier frame of absolutely continuous measures. We let $\mu$ be a compactly supported absolutely continuous probability measure on ${\Bbb R}^n$, so that one can write $d\mu(x) = \varphi(x) dx$ where $\varphi$ is a compactly supported function in $L^1({\Bbb R}^n,dx)$ and the support is $K$. Note that if $\varphi$ satisfies $0<m\leq\varphi(x)\leq M<\infty$ almost everywhere on $K$, then by choosing $R>0$ so that $[-R,R)^n$ contains the support and  $\Lambda$ is a frame spectrum of $L^2([-R,R)^n,dx)$, we can easily check that for any $f\in L^2(K,d\mu)$, we have
$$
m\int|f|^2\varphi dx\leq \sum_{\lambda\in\Lambda}|\int f(x)e^{2\pi i \langle\lambda,x\rangle}\varphi(x)dx|^2\leq M\int|f|^2\varphi dx.
$$
Hence if there exists positive constants  $m,M$ such that  $m\leq\varphi(x)\leq M$ almost everywhere on the support of $\varphi$, then $L^2(K,d\mu)$ admits Fourier frame. Our main result is to obtain the converse.

\medskip

\begin{theorem}\label{th0.1}
  Let  $\mu$ be a compactly supported absolutely continuous probability measure on ${\Bbb R}^n$ with $d\mu(x) = \varphi(x) dx$ and $K$ is its support. Then $L^2(K,d\mu)$ admits a Fourier frame if and only if there exists positive constants  $m,M$ such that  $m\leq\varphi(x)\leq M$ almost everywhere on the support of $\mu$.
\end{theorem}

\bigskip

We can apply the above theorem to characterize those equal weight absolutely continuous self-similar measures which admit a Fourier frame. Let $\{f_j\}_{j=1}^{\ell}$ be an iterated function system with $f_j(x) = \lambda x+d_j$ and $0<\lambda<1$. It is well-known that there exists a unique Borel probability measure $\mu $ satisfying
\begin{equation}\label{eq1.0}
\mu(E) = \sum_{j=1}^{\ell}\frac{1}{\ell}\mu(f_j^{-1}(E))
\end{equation}
for any Borel set $E$. Moreover, the support of the measure is the unique compact set $K$ satisfying $K = \bigcup_{j=1}^{\ell}f_j(K)$. When $\lambda=1/\ell$ and $K$ has positive Lebesgue measure, it is easy to see the invariant $\mu$ in (\ref{eq1.0}) is the Lebesgue measure supporting on $K$. In this case, $K$ is a translational tile in ${\Bbb R}^1$ and $K$ is called a \textit{self-similar tile}. Details of the associated tiling theory can be found in [LW]. It is easy to see from Theorem \ref{th0.1} that if $\mu$ is the Lebesgue measure supported on the self-similar tile, then $L^2(K,d\mu)$ will admit a Fourier frame. We will prove that the converse is also true.
\medskip

\begin{theorem}\label{thm0.2-}
Let $\mu$ be the self-similar measure defined in (\ref{eq1.0}) and $\mu$ is absolutely continuous. If $L^2(K,d\mu)$ admits a Fourier frame, then $\lambda  = \frac{1}{\ell}$, $K$ is a self-similar tile  and the density of $\mu$ is $\chi_K$.
\end{theorem}

\medskip

  For the iterated function system consisting only of $f_{1}(x) = \lambda x$  and $f_2(x)=\lambda x+1-\lambda$, then the unique self-similar measure $\nu_{\lambda}$ defined by (\ref{eq1.0}) is called the  \textit{$\lambda$-Bernoulli convolution}. If $\lambda = 1/n$, it reduces to the standard Cantor measures. It is known that for almost all $\frac{1}{2}\leq\lambda<1$, $\nu_{\lambda}$ is absolutely continuous with respect to the Lebesgue measure [So]. We have the following corollary of Theorem \ref{thm0.2-}.

\medskip

\begin{Cor}\label{thm0.2}
Let $\nu_{\lambda}$ be the Bernoulli convolution on ${\Bbb R}^1$ with support denoted by $K$. If $\nu_{\lambda}$ is absolutely continuous with respect to the Lebesgue measure and $\lambda\neq 1/2$,  then $L^2(K,d\nu_{\lambda})$ cannot admit any Fourier frame. In particular this is true for almost all $\lambda\in(1/2,1)$.
\end{Cor}

\medskip

The question of the existence of orthogonal complex exponentials for the Bernoulli convolution has not been settled completely. Hu and Lau determined those contraction (i.e. $\lambda$) of the Bernoulli convolutions for which there are infinitely many such orthogonal sets [HL]. Dutkay, Han and Jorgensen showed that whenever $\lambda>1/2$, there is no
complete orthogonal complex exponentials [DHJ]. It is also conjectured that there are orthonormal complex exponentials if and only if $\lambda =1/2n$ [{\L}aW].

\bigskip

For the organization of the paper, we first recall some basic properties of Fourier frame and the Beurling density in Section 2. We will also prove some general density results that will be used in our proof. We then prove Theorem \ref{th0.1} in Section 3. We then apply the result to prove Theorem \ref{thm0.2-} and Corollary \ref{thm0.2} in Section 4. Finally, we conclude in Section 5 with some remarks and open questions.

\end{section}

\bigskip
\bigskip

\begin{section}{\bf General density results}

%\begin{Def}
%Let $\mu$ be a probability measure on ${\Bbb R}^n$. $\{e^{2\pi i \langle\lambda,\cdot\rangle}\}_{\lambda\in\Lambda}$  is called a Fourier frame of the $L^2(\mu)$ if there exists $A,B>0$ such that for all $f\in L^2(\mu)$. we have
%
%\begin{equation} \label{Fourier frame}
%A\|f\|^2\leq\sum_{\lambda\in\Lambda}|\int f(x)e^{2\pi i \langle\lambda,x\rangle}d\mu(x)|^2\leq B\|f\|^2.
%\end{equation}
%It is called a Bessel's sequence if only the upper bound of (\ref{Fourier frame}) holds.
%\end{Def}

 In this section, we provide some basic properties of Fourier frame, particularly its connection with the Beurling densities. First, It is easy to see that Fourier frame on subset of ${\Bbb R}^n$ has translational invariance property as in the following proposition.

\medskip

\begin{Prop}\label{thm1.4}
 Let $\Omega$ be a bounded set on ${\Bbb R}^n$ with positive Lebesgue measure. If ${\Lambda}$ is a discrete set on ${\Bbb R}^n$ and $t\in{\Bbb R}^n$. Then $\{e^{2\pi i \langle\lambda,\cdot\rangle}\}_{\lambda\in\Lambda}$ is a  Fourier frame of $L^2(\Omega,dx)$ if and only if it is a  Fourier frame of $L^2(\Omega+t,dx)$.
\end{Prop}

\bigskip

  Let $Q_h(x)  = \prod_{i=1}^{n}[x_i-h/2,x_i+h/2)$ be the cube centered at $x$. The \textit{upper} and \textit{lower Beurling density} of a discrete set $\Lambda$ are defined as follows:
$$
D^{+}\Lambda = \limsup_{h\rightarrow\infty}\sup_{x\in{\Bbb R}^n}\frac{\#(\Lambda\cap Q_{h}(x))}{h^n},\ D^{-}\Lambda = \liminf_{h\rightarrow\infty}\inf_{x\in{\Bbb R}^n}\frac{\#(\Lambda\cap Q_{h}(x))}{h^n}.
$$
A set $\Lambda$ is called \textit{separated} if there exists $\delta>0$ such that $\inf_{x,y\in\Lambda}|x-y|\geq \delta$. It is known that $D^{+}\Lambda<\infty$ if and only if $\Lambda$ is a finite union of separated sequence  [Chr, Lemma 7.1.3].

The study of Fourier frame is closely tied with the density of discrete sets (\cite{[Lan]}, [GR]). Landau gave an important necessary condition on the density for the frame spectrum \cite{[Lan]}. There is also a sufficient condition on ${\Bbb R}^1$ guaranteeing that ${\Lambda}$ is a frame spectrum on an interval (\cite{[Chr]}, \cite{[S]}). We summarize them in the following theorem.

\medskip

\begin{theorem}\label{thm1.1}
 (i) If $\{e^{2\pi i \langle\lambda,\cdot\rangle}\}_{\lambda\in\Lambda}$ is a  Fourier frame of $L^2(\Omega,dx)$ with $\Omega\subset{\Bbb R}^n$ of finite Lebesgue measure, then $D^{-}\Lambda\geq {\mathcal L}(\Omega)$, where ${\mathcal L}$ denotes the Lebesgue measure.

 \medskip

 (ii) If ${\Lambda}$ is a set such that $D^{+}\Lambda<\infty$ and $D^{-}\Lambda>{\mathcal L}(I)$, where $I$ is an interval on ${\Bbb R}^1$, then $\{e^{2\pi i \langle\lambda,\cdot\rangle}\}_{\lambda\in\Lambda}$ is a  Fourier frame of $L^2(I,dx)$.
\end{theorem}

\bigskip

   The study of the Bessel's sequence is more tractable than the Fourier frame. We can actually determine the density criterion for $\Lambda$ to be a Bessel's sequence on $L^2({\Omega},dx)$ .

  \medskip

  \begin{Prop}\label{prop2.1}
  If a set $\Lambda$ satisfies $D^{+}\Lambda<\infty$, then $\{e^{2\pi i \langle\lambda,\cdot\rangle}\}_{\lambda\in\Lambda}$ is a Bessel's sequence of $L^2(\Omega,dx)$ for any bounded subsets $\Omega$ in ${\Bbb R}^n$ with positive Lebesgue measure.
\end{Prop}

\medskip

\noindent{\bf Proof.}
The idea of the proof is essentially the Plancherel-Polya inequality. The case for dimension one can be found in [Y, p.79-83]. The higher dimension case is also known in literature (see e.g. [GR]). We give a short proof for completeness.

\medskip

 As $\Omega$ is bounded, we can find some $T>0$
be such that $\Omega\subset[-T,T]^n$. For any $f\in L^2(\Omega,dx)$, we
have $f\in L^2([-T,T]^n,dx)$. Write $F(\xi) = {\widehat{f}}(\xi) =
\int_{\Omega} f(x)e^{-2\pi i \langle\xi,x\rangle}dx =
\int_{[-T,T]^n}f(x)e^{-2\pi i \langle\xi,x\rangle}dx$. We see that
$F$ is an entire function (on ${\Bbb C}^n$) of exponential type in
the following sense.
$$
|F(x+iy)| \leq  A e^{2\pi T (|y_1|+...+|y_n|)},
$$
where $x =(x_1,...,x_n), \ y = (y_1,...,y_n)\in{\Bbb R}^n$. One can apply the one dimensional Plancherel-Polya Theorem iteratively ([Y,p.79]) (See also [St, Lemma 4.11] for a general statement), we obtain
\begin{equation}\label{2.1}
\int_{-\infty}^{\infty}|F(x+iy)|^2dx\leq e^{4\pi T(|y_1|+...+|y_n|)}\int_{-\infty}^{\infty}|F(x)|^2dx.
\end{equation}

\medskip

Note that $D^{+}\Lambda<\infty$ implies that $\Lambda$ is a finite union of separated sequences. This means $\Lambda = \bigcup_{i=1}^{\ell}\Lambda_i$ with $\Lambda_i$ are separated (i.e. $\lambda_n-\lambda_m\geq \delta_i>0$ for all $m,n$). We can apply a similar argument in [Y,p.82] to prove that
$$
\sum_{i=1}^{\ell}\sum_{\lambda\in\Lambda_i}|F(\lambda)|^2\leq C \|F\|^2.
$$
This is equivalent to $\Lambda$ is a Bessel's sequence of exponentials on  $L^2(\Omega,dx)$.
 \qquad $\Box$

\medskip

The converse of the above proposition is also true. Indeed, it holds for more general measures.

\begin{Prop}\label{thm1.2}
Let $\mu$ be a probability measure on ${\Bbb R}^n$ with support $K$. Suppose $\{e^{2\pi i \langle\lambda,\cdot\rangle}\}_{\lambda\in\Lambda}$ is a Bessel's sequence of $L^2(K,d\mu)$, then $D^{+}\Lambda<\infty$.

%In particular, if $\mu$ be a probability measure on ${\Bbb R}^1$, then $\{e^{2\pi i \langle\lambda,\cdot\rangle}\}_{\lambda\in\Lambda}$ is a Bessel's sequence of $L^2(I)$, where $I$ is any bounded interval in ${\Bbb R}^1$.
\end{Prop}

\medskip

\noindent{\bf Proof.}
Let $\widehat{\mu}$ be the Fourier transform of $\mu$. We suppose on the contrary that $D^{+}\Lambda = \infty$. By Lemma 7.1.3 in \cite{[Chr]}, for any $h>0$ and for any $N\in{\Bbb N}$, we can find some cubes $Q_h(x_N)$ such that
\begin{equation}\label{1.2}
\#(\Lambda\cap Q_{h}(x_N))>N.
\end{equation}
As $\mu$ is a probability measure, $\widehat{\mu}(0) = 1$ and $\widehat{\mu}$ is a continuous function. Hence, there exists $\delta>0$, and $\epsilon>0$ such that whenever $|x|<\delta$, $|\widehat{\mu}(x)|>\epsilon$. Take $h = \delta/2$ in the above, then $\lambda\in\Lambda\cap Q_{h}(x_N)$ implies that $|\lambda-x_N|<\delta$. Consider the function $f_{N}(x) = e^{-2\pi i \langle x_N,x\rangle}$, by (\ref{1.2}),
$$
\sum_{\lambda\in\Lambda}|\int f_N(x)e^{2\pi i \langle\lambda,x\rangle}d\mu(x)|^2\geq\sum_{\lambda\in\Lambda\cap Q_{h}(x_N)}|{\widehat{\mu}}(\lambda-x_N)|^2>N\epsilon^2.
$$
The expression tends to infinity as $N$ tends to infinity. This is a
contradiction.
Hence, we must have $D^{+}\Lambda < \infty$.
\qquad $\Box$
% since
%$$
%\sum_{\lambda\in\Lambda}|\int f_N(x)e^{2\pi i \langle\lambda,x\rangle}d\mu(x)|^2\leq B\|f\|^2 = B.
%$$

%The second statement follows from Lemma 7.1.3 ($\Lambda$ must be relatively separated) and Lemma 7.6.2 in \cite{[Chr]}.

\bigskip

  In the next section, we will need the following simple sufficient condition of a Fourier frame. The main idea of proof has its origin in the fundamental paper of Duffin and Schaeffer [DS], a version of the proof in high dimension was due to [DHSW].

\medskip

\begin{Prop}\label{thm1.5}
Let $\Lambda$ be a set on ${\Bbb R}^n$ such that $D^{+}\Lambda<\infty$. Suppose that for any ${\mathbf k}\in{\Bbb Z}^n$, ${\bf k}+[-\frac{1}{2},\frac{1}{2})^n$ contains at least one element $\lambda_{\bf k}$ in $\Lambda$, then there exists $\epsilon>0$ such that $\{e^{2 \pi i \langle\lambda,\cdot\rangle}\}_{\lambda\in\Lambda}$ is a Fourier frame on $L^2(Q_{\epsilon},dx)$, where $Q_{\epsilon} = [-\frac{\epsilon}{2},\frac{\epsilon}{2})^n$
\end{Prop}

\medskip

\noindent{\bf Proof.}
By Proposition \ref{prop2.1}, the upper bound is satisfied for any $\epsilon>0$. It remains to prove the lower bound is satisfied for some sufficiently small $\epsilon>0$. For notational convenience, we only consider when $n=2$ and high dimension case follows from the same method by considering projection.

\medskip

For each $\lambda_{{\bf k}}\in\Lambda$, we write ${\bf k} = (k_1,k_2)$ and $\lambda_{\bf k} = (\lambda_1,\lambda_2)$. Define $\lambda_{{\bf k}}' = (\lambda_1,k_2)$. We first compare $\lambda_{\bf k}'$ and ${\bf k}$, then we deal with $\lambda_{\bf k}$ and $\lambda_{\bf k}'$. For any $f\in L^2(Q_{\epsilon},dx)$, we let $\xi = (\xi_1,\xi_2)$ and
$$
F(\xi) = \int_{Q_{\epsilon}}f(x)e^{-2\pi i \langle\xi,x\rangle}dx.
$$
Clearly $F$ is analytic in both variables $\xi_1,\xi_2$ and
$$
\frac{\partial^{\ell} F}{\partial\xi_1^{\ell}}(\xi) = \int_{Q_{\epsilon}}f(x) (-2\pi i x_1)^{\ell}e^{-2\pi i \langle\xi,x\rangle}dx.
$$
 For each ${\bf k}\in{\Bbb Z}^2$, using the Taylor expansion at $k_1$ and the Cauchy-Schwarz inequality and $|\lambda_{\bf k}'-{\bf k}|\leq \frac{1}{2}$, we have
 $$
 \begin{aligned}
 |F(\lambda_{\bf k}')-F({\bf k})|^2 =& |\sum_{\ell=1}^{\infty}\frac{\frac{\partial^{\ell} F}{\partial\xi_1^{\ell}}({\bf k})}{{\ell}!}(\lambda_1-k_1)^{\ell}|^2 \\
 \leq& \sum_{\ell=1}^{\infty}\frac{|\frac{\partial^{\ell} F}{\partial\xi_1^{\ell}}({\bf k})|^2}{\ell!}\cdot\sum_{\ell=1}^{\infty}\frac{(1/2)^{2\ell}}{{\ell}!}
 =\sum_{\ell=1}^{\infty}\frac{|\frac{\partial^{\ell} F}{\partial\xi_1^{\ell}}({\bf k})|^2}{\ell!}\cdot(e^{1/4}-1).
 \end{aligned}
 $$
We then also note that by the Parseval's identity  and $f\in L^2(Q_{\epsilon},dx)$,
$$
\sum_{{\bf k}\in{\Bbb Z}^2}|\frac{\partial^{\ell} F}{\partial\xi_1^{\ell}}({\bf k})|^2 = \sum_{{\bf k}\in{\Bbb Z}^2}|\langle (-2\pi i x_1)^{\ell}f, e^{2\pi i \langle {\bf k},\cdot\rangle}\rangle|^2 = \|(-2\pi i x_1)^{\ell}f\|^2\leq (2\pi\epsilon)^{2{\ell}}\|f\|^2.
$$
This shows that
\begin{equation}\label{eq1.4}
 \sum_{{\bf k}\in{\Bbb Z}^2}|F(\lambda_{\bf k}')-F({\bf k})|^2 \leq (e^{1/4}-1)(e^{4\pi^2 \epsilon^2}-1)\|f\|^2.
\end{equation}
Take $\epsilon$ so small that the above constant on the right is small than 1/2.
By the Minkowski's inequality, we have
\begin{equation}\label{eq2.4}
\begin{aligned}
(\sum_{{\bf k}}|F(\lambda_{\bf k}')|^2)^{1/2}\geq (\sum_{{\bf k}}|F({\bf k})|^2)^{1/2}-(\sum_{{\bf k}}|F(\lambda_{\bf k}')-F({\bf k})|^2)^{1/2}
\geq \frac{1}{2}\|f\|.
\end{aligned}
\end{equation}

\bigskip

 Finally, by the Minkowski's inequality and (\ref{eq1.4}), the $\{\lambda_{\bf k}'\}_{{\bf k}}$ is a Bessel's sequence with bound $B = (1+((e^{1/4}-1)(e^{4\pi^2 \epsilon^2}-1))^{1/2})^2$. Hence, repeating the above argument with Taylor expansion with respect to $\lambda_{\bf k}'$, using the Minkowski's inequality and (\ref{eq2.4}) shows that
$$
\begin{aligned}
(\sum_{{\bf k}}|F(\lambda_{\bf k})|^2)^{1/2}\geq& \sum_{{\bf k}}|F(\lambda_{\bf k}')|^2)^{1/2}-(\sum_{{\bf k}}|F(\lambda_{\bf k})-F(\lambda_{\bf k}')|^2)^{1/2}\\ \geq& (\frac{1}{2}- B^{1/2}(e^{1/4}-1)(e^{4\pi \epsilon^2}-1))\|f\|.\\
\end{aligned}
$$
We then choose $\epsilon >0$ even smaller to make the above constant positive, this shows that these ${\lambda}_{\bf k}$'s is a frame spectrum on some small cubes and hence the proof is completed.
 \qquad $\Box$

\end{section}

\bigskip

\bigskip

\begin{section}{\bf Proof of Theorem \ref{th0.1}}

  As indicated in the introduction, we only need to prove the necessity part of Theorem \ref{th0.1}. We will first see that in order for a compactly supported absolutely continuous measure $\mu$ on ${\Bbb R}^n$ to admit a Fourier frame, its density must be bounded below almost everywhere.

\bigskip

%\begin{theorem}\label{thm2.1}
%Let $\mu$ be a compactly supported probability measure on ${\Bbb
%R}^n$ with $L^1$ density $\phi$. Suppose that $L^2(d\mu)$ admits
%Fourier frames, then there exists $m>0$ such that $\phi\geq m$
%almost everywhere on the support of $\mu$.
%\end{theorem}
%
%\medskip

\noindent{\bf Proof of the lower bound.}
Let $d\mu =\varphi(x)dx$, $K = \mbox{supp}\mu$ and $K\subset [-R,R]^n$ for some $R>0$.  We also denote $E_0 = \{x\in K : \varphi(x)\geq 1\}$ and $E_k = \{x\in K : \frac{1}{k+1}<\varphi(x)\leq\frac{1}{k} \}$ for $k\geq 1$, so that
  $$
  K = \bigcup_{k=0}^{\infty}E_k.
  $$
 Suppose that $\varphi$ does not have a lower bound on its support, then $E_k$ has positive Lebesgue measure for infinitely many $k$. By passing to subsequence if necessary, we may assume ${\mathcal L}(E_k)>0$ for all $k$.

\medskip

By assumption, $L^2(K,d\mu)$ has a Fourier frame $\{e^{2\pi i \langle\lambda
,\cdot\rangle}\}_{\lambda\in\Lambda}$, then $D^{+}\Lambda<\infty$ and  hence it is a
Bessel's sequence of $L^2([-R,R]^n,dx)$ by Proposition  \ref{thm1.2}  and \ref{prop2.1}. We now define $f_k =
\chi_{E_{k}}$, note that $E_{k} \subset [-R,R]^n$ and
\begin{equation}\label{eq2.2}
\int_{[-R,R]^n}|f_k\varphi|^2dx = \int_{E_{k}}|\varphi|^2dx\leq \frac{{\mathcal L}(E_{k})}{k^2}\leq \frac{(2R)^n}{k^2}<\infty,
\end{equation}
 Thus $f_k\varphi\in L^2([-R,R]^n, dx)$. Using the Bessel's sequence assumption in
$L^2([-R,R]^n,dx)$ and the Fourier frame lower bound assumption in
$L^2(K,d\mu)$, we obtain
$$
B\int_{[-R,R]^n}|f_k(x)\varphi(x)|^2dx\geq\sum_{\lambda\in\Lambda}|\int_K f_k(x)e^{2\pi i \lambda x}d\mu(x)|^2\geq A\int_{K}|f_k(x)|^2\varphi(x)dx.
$$
Using (\ref{eq2.2}), we find that
$$
\frac{B{\mathcal L}(E_{k})}{k^2}\geq A\int_{K}|f_k(x)|^2\varphi(x)dx = A\int_{E_{k}}\varphi(x)dx\geq \frac{A{\mathcal L}(E_{k})}{k+1}.
$$
This implies that for all $k>0$, $\frac{k+1}{k^2}\geq \frac{A}{B}>0$, which is
a contradiction. Hence, $\varphi$ must be lower bounded almost everywhere. \qquad $\Box$

\bigskip

For the upper bound in the necessity of Theorem \ref{th0.1}, we need to prove several lemmas to compare the Fourier frames of $L^2(K,d\mu)$ and $L^2(E,dx)$ with $E$ is a subset of $K$. In the following, we will use $\varphi|_E$ to denote the restriction of $\varphi$ on $E$ and  $L^{\infty}(E,dx)$ to denote the set of functions that is bounded above almost everywhere on $E$ with respect to the Lebesgue measure.
\medskip

\begin{Lem}\label{lem3.1}
Suppose that $\{e^{2\pi i \langle\lambda,\cdot\rangle}\}_{\lambda\in\Lambda}$ is a
Fourier frame of $L^2(K,d\mu)$, where $d\mu = \varphi dx$ and $K$ is the support of $\mu$.  Then

\medskip

(i) If $E\subset K$ is a set of positive measure and $\varphi|_E\in
L^{\infty}(E,dx)$, then $\{e^{2\pi i \langle\lambda,\cdot\rangle}\}_{\lambda\in\Lambda}$ is a
Fourier frame of $L^2(E,dx)$.

\medskip

(ii) If $F\subset K$ is a set of positive measure such that $\varphi|_F\not\in L^{\infty}(F,dx)$, then $\{e^{2\pi i
\langle\lambda,\cdot\rangle}\}_{\lambda\in\Lambda}$ cannot be a Fourier frame of
$L^2(F,dx)$.

\end{Lem}

\medskip

\noindent{\bf Proof.} First, from the above, we know there exists
$m>0$ such that $\varphi\geq m$ almost everywhere on its support .

\medskip

(i) Let
$f\in L^2(E,dx)$, then we have
$\int|\frac{f(x)}{\varphi(x)}|^2\varphi(x)dx\leq
\frac{1}{m}\int_{E}|f|^2<\infty$. Hence,
$$
\begin{aligned}
\sum_{\lambda\in\Lambda}|\int_{E}f(x)e^{2\pi i \lambda x}dx|^2 = & \sum_{\lambda\in\Lambda}|\int_{E}\frac{f(x)}{\varphi(x)}e^{2\pi i \lambda x}\varphi(x)dx|^2\\
\leq& B\int_{E}|\frac{f(x)}{\varphi(x)}|^2\varphi(x)dx\leq \frac{B}{m}\int_{E}|f(x)|^2dx.
\end{aligned}
$$
This establishes the upper frame bound. For the lower bound, as we have $\varphi\leq M$ almost everywhere on $E$, we have
$$
 \sum_{\lambda\in\Lambda}|\int_{E}\frac{f(x)}{\varphi(x)}e^{2\pi i \lambda x}\varphi(x)dx|^2\geq A \int_{E}|\frac{f(x)}{\varphi(x)}|^2\varphi(x)dx\geq \frac{A}{M} \int_{E}|f(x)|^2dx.
$$

\bigskip

(ii) Suppose $\{e^{2\pi i
\langle\lambda,\cdot\rangle}\}_{\lambda\in\Lambda}$ is a Fourier
frame of $L^2(F,dx)$, we define $D_k = \{x\in K :
k\leq\varphi(x)\leq k+1\}\cap F$ and $f_k = \chi_{D_k}$. As $\varphi|_F\not\in L^{\infty}(F,dx)$,
${\mathcal{L}}(D_k)>0$ for infinitely many $k$. We may assume that
it holds for all $k$. Note that
$$
\int|f_k(x)\varphi(x)|^2dx = \int_{D_k}|\varphi(x)|^2dx\leq (k+1)^2{\mathcal{L}}(D_k)<\infty.
$$
Hence, by the Fourier frame assumption on $L^2(F,dx)$ and the Bessel's sequence assumption  on $L^2(K,d\mu)$, we obtain
$$
A k^2{\mathcal{L}}(D_k) \leq A\int|f_k(x)\varphi(x)|^2dx\leq B\int|f_k(x)|^2\varphi(x)dx  \leq B (k+1) {\mathcal{L}}(D_k).
$$
This implies that $\frac{k+1}{k^2}\geq \frac{A}{B}$ for all $k$ which is a contradiction.
\qquad $\Box$
%As $\{e^{2\pi i \lambda\cdot}\}_{\lambda\in\Lambda}$ is a Fourier
%frame of $L^2(E_N)$, by the Landau density theorem (Theorem \ref{thm1.1}),
%we have $D^{-}\Lambda\geq {\mathcal L}(E_N)$. As $E_N$ are
%increasing sequence of sets and $\bigcup_{N}E_N = K$, the inequality
%follows by $\lim_{N\rightarrow\infty}{\mathcal L}(E_N) = {\mathcal
%L}(K)$.
\medskip

\begin{Lem}\label{lem3.2}
Let $\{e^{2\pi i \langle\lambda,\cdot\rangle}\}_{\lambda\in\Lambda}$ be a
Fourier frame of $L^2(K,d\mu)$, where $d\mu = \varphi dx$ and $K =$supp$\mu$. Suppose that $\varphi\not\in L^{\infty}(K,dx)$, then $\{e^{2\pi i \langle\lambda,\cdot\rangle}\}_{\lambda\in\Lambda}$ cannot be a Fourier frame of $L^2(Q,dx)$ on any cube $Q$ in ${\Bbb R}^n$.
\end{Lem}

\medskip

\noindent{\bf Proof.}
  For ${\bf k} = (k_1,...,k_n)$, let $I_{{\bf k},r}$ be the dyadic cube $[\frac{k_1}{2^r},\frac{k_1+1}{2^r})\times...\times[\frac{k_n}{2^r},\frac{k_n+1}{2^r})$. Then $\{I_{{\bf k},r}: {\bf k}\in{\Bbb Z}^n, r\in{\Bbb Z}\}$ is the set of all dyadic cubes in ${\Bbb R}^n$.

 To prove the statement, it suffices to prove that $\{e^{2\pi i \langle\lambda,\cdot\rangle}\}_{\lambda\in\Lambda}$ cannot be a Fourier frame on any dyadic cubes . By Proposition \ref{thm1.4}, it suffices to show that for each integer $r$, $\{e^{2\pi i \langle\lambda,\cdot\rangle}\}_{\lambda\in\Lambda}$ cannot be a frame on at least one dyadic cube with side length $2^{-r}$. Let $r$ be given, we note that $K$ is compact,  and thus $K$ is covered by a finite number of dyadic cubes $I_{{\bf k}_i,r}$ of length $2^{-r}$.
As $\varphi\not\in L^{\infty}(K,dx)$, there exists some dyadic cubes $ Q = I_{{\bf k}_i,r}$  such that $\varphi|_{Q\cap K}\not\in L^{\infty}(Q\cap K,dx)$ and ${\mathcal L}(Q\cap K)>0$.  By Lemma \ref{lem3.1}(ii), $\{e^{2\pi i \langle\lambda,\cdot\rangle}\}_{\lambda\in\Lambda}$ cannot be a Fourier frame on $Q\cap K$. This means $\{e^{2\pi i \langle\lambda,\cdot\rangle}\}_{\lambda\in\Lambda}$ cannot be a Fourier frame on $Q$ since ${\mathcal L}(Q\cap K)>0$.
 \qquad $\Box$

\bigskip

Combining the above lemmas and the density results in Theorem \ref{thm1.1}, we can now prove the existence of the upper bound.

%\begin{theorem}\label{thm3.1}
%Let $\mu$ be a compactly supported probability measure on ${\Bbb
%R}^n$ with $L^1$ density $\phi$. Suppose that $L^2(d\mu)$ admits
%Fourier frames, then there exists $M>0$ such that $\phi\leq M$
%almost everywhere on the support of $\mu$.
%\end{theorem}

\medskip

\noindent{\bf Proof of the upper bound.} We argue by contradiction. Suppose that there does not exist $M>0$ such that  $\varphi\leq M$ almost everywhere on the support $K$  and  $L^2(K,d\mu)$ still admits a
Fourier frame $\{e^{2\pi i \langle\lambda,\cdot\rangle}\}_{\lambda\in\Lambda}$. We let
$$
E_N = \{x\in K: \varphi(x)\leq N\}.
$$
then by Lemma \ref{lem3.1}(i),  $\{e^{2\pi i \langle\lambda,\cdot\rangle}\}_{\lambda\in\Lambda}$ is a Fourier
frame of $L^2(E_N,dx)$. By the Landau's density theorem (Theorem \ref{thm1.1}),
we have $D^{-}\Lambda\geq {\mathcal L}(E_N)$. As $E_N$ are
increasing sequence of sets and $\bigcup_{N}E_N = K$, we have
$$
D^{-}\Lambda\geq {\mathcal L}(K).
$$
As $d\mu$ is absolutely continuous, so the support must have positive Lebesgue measure and hence $D^{-}\Lambda\geq c>0$.

\medskip

As $D^{-}\Lambda\geq c$, we can find large $L$ so that for any $x\in{\Bbb R}^n$, $x+[-\frac{L}{2},\frac{L}{2})^n$ contains at least one point of $\Lambda$. Define $\Gamma = L{\Bbb Z}^n$, Then for any $\gamma\in\Gamma$, there exists $\lambda_{\gamma}\in\Lambda$ such that
\begin{equation}\label{eq3.2}
\lambda_{\gamma}\in \gamma+[-\frac{L}{2},\frac{L}{2})^n.
\end{equation}
Denote $\Lambda' = \{\lambda_{\gamma}:\gamma\in\Gamma\}$, we have $\Lambda'\subset\Lambda$ so that $D^{+}\Lambda'<\infty$. Moreover, by (\ref{eq3.2}) and the definition of $\Gamma$, every cube ${\bf k}+[\frac{1}{2},\frac{1}{2})^n$ with $k\in{\Bbb Z}^n$ has one point in $\frac{1}{L}\Lambda'$. By Proposition \ref{thm1.5}, $\{e^{2\pi i \langle\frac{1}{L}\lambda',\cdot\rangle}\}_{\lambda'\in\Lambda'}$ is a Fourier frame of $L^2(Q_{\epsilon},dx)$, where $Q_{\epsilon} = [-\frac{\epsilon}{2},\frac{\epsilon}{2})^n$ and $\epsilon$ is sufficiently small. This implies that  $\Lambda$ will generate a Fourier frame of the $L^2$ space of a cube of side length $\frac{\epsilon}{L}$. This is a contradiction to Lemma \ref{lem3.2}. Thus, we conclude that $\varphi$ must be bounded above almost everywhere. This completes the proof of the upper bound and hence Theorem \ref{th0.1}.
 \qquad $\Box$
\end{section}

\bigskip
\bigskip

\begin{section}{\bf Self-similar measures}

In this section, we consider the iterated function system $f_j(x) = \lambda x+d_j$, for $j=1,...,\ell$ and $0<\lambda<1$. Let ${\mathcal D} = \{d_j:j=1,...,\ell\}$, it is well-known that there exists a unique Borel probability measure $\mu = \mu(\lambda,\ell,{\mathcal D})$ satisfying
\begin{equation}\label{eq4.0}
\mu(E) = \sum_{j=1}^{\ell}\frac{1}{\ell}\mu(f_j^{-1}(E))
\end{equation}
for any Borel set $E$. Moreover, the support of this measure is the unique compact set $K$ satisfying $K = \bigcup_{j=1}^{\ell}f_j(K)$. Explicitly, we can write
\begin{equation}\label{eq4.1}
K = \{\sum_{j=0}^{\infty}\lambda^jd_j: d_j\in{\mathcal D}\}.
\end{equation}
By a suitable translation, we can assume $0= d_1<d_2<...<d_{\ell}$ so that $0$ is in the support $K$ from (\ref{eq4.1}). There are literatures determining whether such measures are absolute continuous (see e.g., [DFW], [LLR]). In Theorem \ref{thm4.1}, we characterize this kind of absolutely continuous measures which admits a Fourier frame.  We start with a lemma.

\medskip

\begin{Lem}\label{lem4.1}
Let $\mu = \mu(\lambda,\ell,{\mathcal D})$ be the self-similar measure defined in (\ref{eq4.0}) and $\mu$ is absolutely continuous with respect to the Lebesgue measure. Then $\lambda\geq1/\ell$ and
\begin{equation}\label{eq4.0+}
\liminf_{n\rightarrow\infty}\frac{{\mathcal L}(K\cap[0,\lambda^n))}{\lambda^n}\geq c>0
\end{equation}
for some constant $c$.
\end{Lem}

\medskip

\noindent{\bf Proof.}
Since $\mu$ is absolutely continuous, ${\mathcal L}(K)>0$. By taking the Lebesgue measure to $K = \bigcup_{j=1}^{\ell}f_i(K)$, we have ${\mathcal L}(K)\leq \lambda\ell{\mathcal L}(K)$. Hence, $\lambda\geq1/\ell$ follows.

\medskip

To prove (\ref{eq4.0+}), we know from (\ref{eq4.1}) that $K$ lies in the non-negative real line, so we can take $N$ so that $\lambda^N K\subset [0,1)$. Hence, $\lambda^{N+n}K\subset [0,\lambda^n)$. By iterating the system $N+n$ times and noting that $d_1=0$, it is easy to see that $K\supset\lambda^{N+n}K$. Thus,
$$
K\cap[0,\lambda^n)\supset\lambda^{N+n}K\cap[0,\lambda^n) = \lambda^{N+n}K.
$$
Taking the Lebesgue measure, we have ${\mathcal L}(K\cap[0,\lambda^n))\geq \lambda^{N+n}{\mathcal L}(K)$. (\ref{eq4.0+}) follows by letting $c = \lambda^{N}{\mathcal L}(K)$ ($>0$).
 \qquad $\Box$

\bigskip

We now state the main theorem in this section.

\begin{theorem}\label{thm4.1}
Let $\mu = \mu(\lambda,\ell,{\mathcal D})$ be the self-similar measure defined in (\ref{eq4.0}) and is absolutely continuous with supp$\mu=K$, the self-similar set. If $L^2(K,d\mu)$ admits a Fourier frame, then $\lambda  = \frac{1}{\ell}$, the density of $\mu$ is $\chi_K$ and $K$ is a self-similar tile.
\end{theorem}

\medskip

\noindent{\bf Proof.}
We first consider $\mu[0,\lambda^n)$. By applying (\ref{eq4.0}), we have
\begin{equation}\label{eq4.2}
\mu[0,\lambda^n) = \frac{1}{\ell}\mu[0,\lambda^{n-1})+\sum_{j=2}^{\ell}\frac{1}{\ell}\mu[-\frac{d_j}{\lambda},\lambda^{n-1}-\frac{d_j}{\lambda}).
\end{equation}
Taking $N$ large enough so that for all $n\geq N$, we have $\lambda^{n-1}-\frac{d_j}{\lambda}<0$ for all $j=2,..., \ell$. By $(\ref{eq4.2})$ and noting that the support of $\nu_{\lambda}$ lies in the non-negative real line, we conclude that for all $n\geq N$,
\begin{equation}\label{eq4.3}
\mu[0,\lambda^n)= \frac{1}{\ell}\mu([0,\lambda^{n-1})) = \cdots =\frac{1}{{\ell}^{n-N}}\mu[0,\lambda^{N})\leq C\cdot(\frac{1}{\ell})^n,
\end{equation}
where $C$ is independent of $n$.

\medskip

 Let $\varphi$ be the density of $\mu$ . As $\mu$ admits a Fourier frame, by Theorem \ref{th0.1}, we have $\varphi\geq m>0$ almost everywhere on its support. By (\ref{eq4.3}) and Lemma \ref{lem4.1}, for $n$ large
 $$
 C(\frac{1}{\ell})^n\geq \int_{0}^{{\lambda}^n}\varphi(x)dx\geq m{\mathcal L}(K\cap[0,\lambda^n))\geq mc \lambda^n.
 $$
  This means that $\lambda\leq 1/\ell$. Combining with Lemma \ref{lem4.1}, $\lambda = 1/\ell$. As ${\mathcal L}(K)>0$, $K$ must be a self-similar tile on ${\Bbb R}^1$ and the density is clearly $\chi_K$.
 \qquad $\Box$

 \bigskip

 As a corollary, we consider the $\lambda$-Bernoulli convolution $\nu_{\lambda}$, which is  the unique self-similar measure defined by the iterated function system $f_{1}(x) = \lambda x$  and $f_2(x)=\lambda x+1-\lambda$ as in (\ref{eq4.0}).

\medskip

\begin{Cor}\label{thm2.2}
Let $\nu_{\lambda}$ be the Bernoulli convolution. If $\nu_{\lambda}$ is absolutely continuous with respect to the Lebesgue measure,  then $L^2(K,d\nu_{\lambda})$ cannot admit any Fourier frame if $\lambda\neq1/2$. In particular, for almost all $\lambda\in(1/2,1)$, $L^2(K,d\nu_{\lambda})$ cannot admit any Fourier frame.
\end{Cor}

\medskip

\noindent{\bf Proof.}
By Theorem \ref{thm4.1}, if $L^2(K,d\nu_{\lambda})$ admits a Fourier frame, then $\lambda = 1/2$. This shows the first statement. The second statement is a direct consequence of the fact that for almost all $\lambda\geq1/2$,  $\nu_{\lambda}$ is absolutely continuous [So].
\qquad $\Box$
\end{section}

\bigskip
\bigskip

\begin{section}{\bf Remarks and open questions}

Let $\mu$ be a Borel probability measure on ${\Bbb R}^n$ and let $\widehat{\mu}$ be the Fourier transform of $\mu$, then by putting $e^{2\pi i \langle \xi,\cdot\rangle}$ into the definition of the Fourier frame, we obtain a necessary condition for the existence of a Fourier frame.
\begin{equation} \label{Fourier frame condition}
A\leq\sum_{\lambda\in\Lambda}|\widehat{\mu}(\xi+\lambda)|^2\leq B,
\end{equation}

\medskip

When the frame is an orthonormal basis, then $A = B= 1$ in (\ref{Fourier frame condition}). It is known that this identity is sufficient for completeness of the orthogonal set  $\{e^{2\pi i \langle\lambda,\cdot\rangle}\}_{\lambda\in\Lambda}$ [JP]. When we only assume the frame condition, it was asked by Dutkay and Jorgensen [DJ] that whether the inequality is still sufficient. In the following example, we show that the answer is negative.

\medskip

\begin{Example}
Let $m$ be the Lebesgue measure supported on $[-1/2,1/2]$, and $\mu = m\ast m$. Then $\mu$ does not admit any Fourier frame but $0<A\leq\sum_{n\in{\Bbb Z}}|\widehat{\mu}(x+n)|^2\leq 1$.
\end{Example}

\medskip

\noindent{\bf Proof.}
By a direct calculation, we see that  $m\ast m$ is absolutely continuous with density $f(x )  = 1-|x|$. Thus, $f$ is not  bounded below almost everywhere. This means $m\ast m$ does not admit Fourier frame by Theorem \ref{th0.1}.

\medskip

On the other hand,
$$
\sum_{n\in{\Bbb Z}}|\widehat{\mu}(x+n)|^2 = \sum_{n\in{\Bbb Z}}|\widehat{m}(x+n)|^4\leq \sum_{n\in{\Bbb Z}}|\widehat{m}(x+n)|^2 =1.
$$
We also note that $\sum_{n\in{\Bbb Z}}|\widehat{\mu}(x+n)|^2$ is integer periodic, so
$$
\sum_{n\in{\Bbb Z}}|\widehat{m}(x+n)|^4\geq \inf_{x\in(-1/2,1/2]}\sum_{n\in{\Bbb Z}}|\widehat{m}(x+n)|^4 \geq  \inf_{x\in(-1/2,1/2]}|\widehat{m}(x)|^4= (\frac{\sin(\pi/2)}{\pi/2})^4 >0.
$$
This completes the proof.
\qquad $\Box$

\bigskip

There are many questions remaining open. Another closely related objects of Fourier frame is the \textit{exponential Riesz basis}. Recall that $\{e_n\}$ is a \textit{Riesz basis} in a separable Hilbert space $H$ if it is complete and for any $\{c_n\}\in\ell^2$, there exists $f\in H$ such that
$$
\langle f,e_n\rangle =c_n.
$$
A Riesz basis is also equivalent to an \textit{exact frame}, i.e. a frame that fails to be a frame if one of the vectors is removed. It is also clear that an orthonormal basis is a Riesz basis.  An exponential Riesz basis is the Riesz basis of the form $e^{2\pi i \langle\lambda,\cdot\rangle}$ in $L^2(K,d\mu)$, for some compactly supported probability measure $\mu$. In view of Theorem \ref{th0.1}, we ask

 \bigskip

 \noindent\textbf{Q1:}  Can we classify the densities $\varphi$ that $L^2(K,\varphi dx)$ admits an exponential Riesz basis or exponential orthonormal basis?

 \bigskip

  The question for the exponential orthonormal basis is a generalization of the Fuglede conjecture [Fu]: $L^2(\Omega,dx)$ admits an exponential orthonormal basis if and only if $\Omega$ is a translational tile. Although it was proved to be false in general [T], the exact relationship with tiling is still widely open.

\bigskip

  In section 4, we have shown that the only equal weight absolutely continuous self-similar measure admitting a Fourier frame is the self-similar tiles. Up to now, except the case of self-affine tiles, we cannot find any example of self-similar (or self-affine) measure on ${\Bbb R}^n$ that is absolutely continuous with density bounded above and below almost everywhere on its support. It is conjectured that Theorem \ref{thm4.1} will hold more generally:

\bigskip

 \noindent\textbf{Q2:}  Is it true that the only absolutely continuous self-affine measures on ${\Bbb R}^n$ admitting a Fourier frame is the characteristic function of self-affine tiles?

 \bigskip

This paper focuses only on the absolutely continuous measures, the question becomes more difficult when the measure is singular. One reason is that there is a lack of Fourier duality theory for general singular measures, it is then hard to produce a good necessary condition in terms of Beurling densities as in Theorem \ref{thm1.1}. Recently, Dutkay \textit{et al} [DHSW] found some necessary conditions for the existence of Fourier frame of fractal measures with the open set condition in terms of the \textit{Beurling dimension}. They also showed that all fractal measures arising from iterated function system with equal contraction admit Bessel's exponential sequence of some positive Beurling dimension [DHW]. Despite such intensive studies, there is up to now no examples of fractal self-affine measure admitting Fourier frame but not exponential orthonormal basis. Here, we post the following question:

\bigskip

 \noindent\textbf{Q3:}  Can we classify the singular measures which admits a Fourier frame?

\bigskip

\noindent{\bf Acknowledgement.}  I would like thank Professor Ka-Sing Lau for his guidance, comments and encouragement when preparing for this paper. I would also like to thank Professor De-Jun Feng for providing the idea of proving Theorem \ref{thm4.1}, which greatly improved the result when compared to the original version.
\end{section}


\bigskip
\bigskip


\begin{thebibliography}{9999}

\bibitem [Chr] {[Chr]}
{\sc O. christensen}, {\it An Introduction to Frames and Riesz Bases}, Applied and Numerical Harmonic Analysis. Birkh\"{a}user Boston Inc., Boston, MA, 2003.

\bibitem [DFW] {[DFW]}
{\sc X.R. Dai, D.J. Feng and Yang Wang}, {\it Refinable functions
with non-integer dilations.}, J. Func. Anal., 250 (2007), 1-20.


\bibitem [DHJ]{[DHJ]}
{\sc D. Dutkay, D.G. Han, and P. Jorgensen}, {\it Orthogonal
exponentials, translations and Bohr completions}, J. Funct. Anal.,
257 (2009), 2999-3019.



\bibitem [DHSW]{[DHSW]}
{\sc D. Dutkay, D.G. Han, Q.Y. Sun and E. Weber}, {\it On the
Beurling dimension of exponential frames}, Adv in Math., 226 (2011),
 285-297.


\bibitem [DHW]{[DHW]}
{\sc D. Dutkay, D.G. Han, and E. Weber}, {\it Bessel sequence of exponential on fractal measures}, preprint.

\bibitem [DJ]{[DJ]}
{\sc D. Dutkay, and P. Jorgensen}, {\it Affine fractals as boundaries and their harmonic analysis}, preprint.


\bibitem [Fu]{[Fu]}
{\sc B. Fuglede}, {\it Commuting self-adjoint partial differential operators and a group theoretic problem}, J. Funct. Anal., 16 (1974), 101-121.

\bibitem [G] {[G]}
{\sc K. Gr\"{o}chenig}, {\it Foundations of time-frequency analysis}, Applied and Numerical Harmonic Analysis. Birkh\"{a}user, Boston, Basel, Berlin, 2001.

\bibitem [GR] {[GR]}
{\sc K. Gr\"{o}chenig and H. Razafinjatovo}, {\it On Landau's necessary density conditions for sampling and interpolation of band-limited functions}, J. Lonon Math. Soc. 54 (1996), 557-565.


\bibitem [HL] {[HL]}
{\sc T.Y. Hu and K.S. Lau}, {\it Spectral property of the Bernoulli convolutions.}, Adv. Math., 219 (2008), 554-567.



\bibitem [Ja] {[Ja]}
{\sc S. Jaffard}, {\it A density theorem for frames of complex exponentials}, Michigian Math. J., 38 (1991), 339-348.

\bibitem [JP] {[JP]}
{\sc P. Jorgensen and S. Pedersen}, {\it Dense analytic subspaces in fractal $L^2$ spaces.}, J. Anal. Math., 75 (1998), 185-228.



\bibitem [Lan] {[Lan]}
{\sc H. Landau}, {\it Necessary Density Conditions for Sampling and
Interpolation of Certain Entrie Functions}, Acta Math., 117 (1967),
37-52.

\bibitem [{\L}aW] {[LaW]}
{\sc I. {\L}aba and Y. Wang}, {\it On spectral Cantor measures}, J.
Funct. Anal., 193 (2002), 409 - 420.

\bibitem [LLR] {[LLR]}
{\sc C.K. Lai, K.S. Lau and H. Rao}, {\it  Spectral structure of digit sets of self-similar tiles on ${\Bbb R}^1$}, preprint.

\bibitem [LW] {[LW]}
{\sc J. Lagarias and Y. Wang}, {\it Self-Affine tiles in ${\Bbb R}^n$}, Adv. in Math., 121 (1996), 21 - 49.

\bibitem [OS] {[OS]}
{\sc J. Ortega-Cerd\`{a} and K. Seip}, {\it Fourier frames}, Ann of Math., 155 (2002),
789-806.

\bibitem [S] {[S]}
{\sc K. Seip}, {\it On the Connection between Exponential Bases and
Certain Related Sequences in $L^2(-\pi,\pi)$}, J. Funct, Anal., 130
(1995), 131 -- 160.


\bibitem [So] {[So]}
{\sc B. Solomyak}, {\it On the Random Series $\sum\pm\lambda^n$ (an Erd\"{o}s problem)}, Ann. of Math., 142 (1995), 611-625.

\bibitem [St] {[St]}
{\sc E. Stein}, {\it Introduction to Fourier Analysis on Euclidean Spaces}, Princeton University Press, 1971.

 \bibitem [T] {[T]} {\sc T. Tao}, {\it Fuglede's conjecture is false in 5 or higher dimensions}, Math. Res. Letter, 11 (2004), 251-258.

\bibitem [Y] {[Y]}
{\sc R. Young}, {\it An Introduction to Nonharmonic Fourier Series}, San Diego : Academic Press, Rev. 1st ed, c2001.

\end{thebibliography}
\end{document}